\documentclass[12pt,leqno]{article}
\usepackage[francais]{babel} 
\usepackage[latin1]{inputenc}
\usepackage[T1]{fontenc}
\usepackage{amsmath}
\usepackage{amsfonts}
\usepackage{amssymb}

\newcommand\N{\mathbb{N}}
\newcommand\R{\mathbb{R}}

\newcommand\CD{\mathcal{D}}
\newcommand\CN{\mathcal{N}}
\newcommand\CM{\mathcal{M}}
\newcommand\CP{\mathcal{P}}
\newcommand\CF{\mathcal{F}}

\newcommand\rg{\rightarrow}

\begin{document}

\title{Conditions pour que les entiers de\\
Beurling aient une densité}

\author{Jean--Pierre Kahane}

\date{23 Octobre 2016}

\maketitle

\textbf{Résumé :}
En 1977 H.G.Diamond donna une condition portant sur les nombres premiers généralisés de Beurling qui entraîne que les entiers correspondants aient une densité. Nous donnons une nouvelle preuve que cette condition est suffisante (Théorème 1) et nous montrons qu'elle n'est pas nécessaire (Théorème 2 et Exemples), mais qu'elle est néanmoins très près d'être nécessaire et suffisante (Théorème 3). Les preuves des Théorèmes 1 et 2 reposent sur l'analyse de Fourier.

\vspace{2mm}

\textbf{Abstract :}   
In 1977 H.G.Diamond gave a condition on Beurling's generalized prime numbers in order that the corresponding generalized integers have a density. We give a new proof of this condition (Théorème~1) and a proof that it is not necessary (Théorème 2 and Examples). However, it is very near to be necessary (Théorème 3). Both proofs of Theorems 1 and 2 rely on Fourier analysis.

\vspace{2mm}

\textbf{Mots clés} :  Key words
Beurling, Diamond, generalized prime undintegers, Fourier methods  in numbertheory

\vskip4mm

\section{\kern-4mm Théorème de Diamond et analyse har\-mo\-nique}

Harold Diamond a donné en 1977 une condition portant sur les nombres premiers de Beurling pour que les entiers correspondants aient une densité. Nous allons la rappeler, la commenter, donner par analyse de Fourier une nouvelle démonstration, et mettre en évidence d'autres conditions. On verra le rôle de l'analyse harmonique dans les théorèmes~1 et 2, et une autre approche pour le théorème~3. Disons d'abord de quoi il s'agit \cite{1,2,3}.

Dans la théorie de Beurling, toute suite croissante et tendant vers l'infini de nombres réels contenus dans la demi--droite ouverte $]1,\infty[$ est considérée comme une suite de \og nombres premiers \fg{,} et les \og entiers \fg{} correspondants sont engendrés par leur multiplication, avec une multiplicité égale au nombre de représentations comme produits de \og premiers \fg{.} On note $\CP$ l'ensemble des \og nombres premiers \fg{} et $\CN$ le multiensemble des \og entiers \fg{,} c'est--à--dire l'ensemble des \og entiers \fg{} comptés avec leur multiplicité, et on désigne
 par $\CP(y)$ et $\CN(x)$ ; les sauts de $\CP(y)$ et de $\CN(x)$ donnent la position et éventuellement la multiplicité des éléments de $\CP$ et de $\CN$. On dira que $\CN$ a une densité (on sous--entendra toujours : finie et non nulle) quand $\CN(x)/x$ a une limite finie et non nulle quand~$x\rg \infty$ :
\begin{equation}
    \CN(x) = Dx+o(x)\ \ (x\rg \infty),\ 0<D<\infty \end{equation} 
et la limite $D$ est la densité. La condition de Diamond est
\begin{equation}
    \int_e^\infty \Big|\CP(y) - \frac{y}{\log y}\Big| \frac{dy}{y^2} < \infty\,.
\end{equation}
Naturellement, on peut remplacer $\frac{y}{\log y}$ par $\CP_i(y)$ quand
$$
    \int_e^\infty \Big|\CP_i(y) - \frac{y}{\log y}\Big| \frac{dy}{y^2} < \infty\,.
$$    
et c'est ce que nous ferons dans deux cas : $\CP_1(y)=\pi(y)$, la fonction de comptage des nombres premiers usuels, et $\CP_o(y)=\tau(y)$ comme on verra plus tard, avec la formule~(5).

Voici deux observations simples

\vspace{2mm}

1. \textbf{Dans le cas où les \og nombres premiers \fg{} de Beurling sont des nombres premiers usuels, la condition de Diamond est nécessaire et suffisante. }

En effet, on obtient $\CP$ en supprimant de $\CP_1$, l'ensemble des nombres premiers usuels, l'ensemble $\CP' =\CP_1-\CP$, et la suppression de chaque $p\in \CP'$ multiplie la densité des entiers par $1-\frac{1}{p}$. La densité finale est $\prod\limits_{p\in \CP'}(1-\frac{1}{p})$, qui est non nulle si et seulement si $\sum\limits_{p\in \CP'} \frac{1}{p}<\infty$, et c'est précisément la condition de Diamond.

\vspace{2mm}

2. \textbf{Dans le cas où $\CP$ est multiplicativement libre et contient $\CP_1$, il en est de même.}

En effet, à chaque stade, on multiplie la densité des entiers par $(1-\frac{1}{p})^{-1}$ $(p\in \CP-\CP_1)$. 

Dans ces deux cas, le signe de $\CP(y)-\pi(y)$ est constant. Nous verrons plus loin, à la suite du théorème~2, des exemples où ce signe varie de telle sorte que la condition~(2) est violée et que cependant $\CN$ a une densité. Nous verrons également, au théorème~3, comment élargir ces observations en un cadre qui rend la condition~(2) nécessaire et suffisante.

\vspace{2mm}

\textbf{L'extension de Beurling}

La relation entre $\CP$ et $\CN$ s'exprime par la formule
$$
\sum_{\nu\in \CN}
 \nu^{-s}= \prod_{p\in \CP}(1-p^{-s})^{-1}
 $$
 qu'on peut aussi écrire
 \begin{equation}
 \int x^{-s} d\,\CN(x) = \exp \int \log (1-y^{-s})^{-1} d\, \CP(y)
 \end{equation}
 Beurling élargit le cadre : pour lui, et pour nous désormais, $\CP(y)$ et $\CN(x)$ sont deux fonctions croissantes sur $]1,\infty[$, avec $\CN(1)=1$, liées par la relation~(3). Quand $\CN(x)$ vérifie (1), la première intégrale dans (3) définit une fonction analytique $Z(s)$ dans le demi--plan ouvert $\sigma>1$ $(s=\sigma+it)$. Nous supposerons toujours 
 \begin{equation}
 \CP(y)=o(y)\quad (y\rg \infty)\,.
 \end{equation}
 Alors la seconde intégrale s'écrit aussi
$$
s\ \int \frac{y^{-s-1}}{1-y^{-s}} \CP(y)dy
$$
qui est une fonction analytique dans le demi--plan $\sigma>1$.

Nous nous attacherons au triplet $\CP(y)$, $\CN(x)$, $Z(s)$. Comme triplet de référence nous prendrons d'abord
\begin{equation}
\CP_o(y)\! =\tau(y)\!\! =\! \int_1^y (1-\xi^{-1})\?log^{-1}\xi\, d\, \xi\,,\, \CN_o(x)\!=\!x\,1_{[1,\infty[}(x)\,,\ Z_o(1)\!=\!\frac{s}{s-1}
\end{equation}
et nous poserons
$$
\CP(y)-\CP_o(y)=a(y)\,.
$$
La condition de Diamond s'écrit
\begin{equation}
\int \frac{|a(y)|}{y^2} d y<\infty
\end{equation}
et l'hypothèse (4) s'écrit $a(y)=o(y)$ $(y\rg \infty)$. Introduisons
\begin{eqnarray}
&A(s) &=\int y^{-s} d a(y) = s \int y^{-s-1}a(y)dy\\
&B(s)&= \int\log(1-y^{-s})^{-1} da(y) =s \int \frac{y^{-s-1}}{1-y^{-s}} a(y) dy\\
&&=A(s) +\frac{1}{2} A(2s) +\frac{1}{3} A(3s)  + \cdots \nonumber\\
&C(s) &=\exp B(s)\,.
\end{eqnarray}
Alors
\begin{equation}
Z(s) =Z_o(s) C(s) = \frac{s}{s-1} C(s)
\end{equation}

Sous la condition (6), les fonctions $A(s)$, $B(s)$ et $C(s)$ sont bien définies et continues dans le demi--plan fermé $\sigma\ge 1$, et $A(s) =O(|s|)$, $B(s) =O(|s|)$ et $C(s)=\exp O(|s|)$ quand $|s|\rg \infty$. De plus, cette condition entraîne une propriété essentielle des trois fonctions $A(1+it)$, $B(1+it)$ et $C(1+it)$ : elles appartiennent localement à l'algèbre de Wiener $W=\CF L^1(\R)$\footnote{Sur l'appartenance locale à $W$ et sur les estimations qui suivent, on peut consulter, entre autres, les pages~11 à 13 du livre \textit{Séries de Fourier absolument convergentes}, de J.--P.~Kahane, Ergebnisse der Mathematik 50, Springer 1970}, avec des normes dans les algèbres $W(I)$ (algèbres des restrictions à des intervalles $I=[\theta,\theta+1]$) respectivement $O(\theta)$ et $\exp O(\theta)$ $(\theta\rg \infty)$. Ainsi, en posant
\begin{equation}
\gamma_\varepsilon(t) = \exp (-\varepsilon^2 t^2/2)
\end{equation}
compte tenu de la décroissance des normes de $\gamma_\varepsilon(\cdot)$ dans les $W(I)$ quand $\theta\rg \infty$, on voit par partition de l'unité que pour tout~$\varepsilon>0$
\begin{equation}
C(1+it)\gamma_\varepsilon(t) \in \CF L^1(\R)\,.
\end{equation}

Pour étudier $\CN(x)$ à partir de $\CP(y)$, on peut songer à la formule classique
 $$
 \CN(x) =\frac{1}{2\pi i} \int_{\sigma- i \infty}^{\sigma+ i\infty} x^sZ(s) \frac{ds}{s}\,,
 $$
 mais la formule (10) nous dit seulement que, sous l'hypothèse (6), $Z(s) =\exp O(|s|)$ $(|s| \rg \infty)$. On est donc amené à introduire $\gamma_\varepsilon(t)$ et écrire
\begin{equation}
\CN_\varepsilon(x)=\frac{1}{2\pi} \int_{-\infty}^\infty  x^{\sigma+it} \frac{Z(\sigma+it)}{\sigma+it} \gamma_\varepsilon (t) dt \quad (\sigma>1)
\end{equation}
l'intégrale étant une intégrale de Lebesgue. Il n'est pas permis d'écrire (13) pour $\sigma=1$, à cause du pôle de $Z_o(s)$ au point~1. Mais~on~a
\begin{equation}
\frac{d}{du}(e^{-\sigma u}\CN_\varepsilon(e^u)) = \frac{1}{2\pi} \int_\R it\, e^{iut} \frac{Z(\sigma+it)}{\sigma+it} \gamma_\varepsilon(t)dt
\end{equation}
qui est valable pour $\sigma\ge 1$, et pour $\sigma=1$ on obtient
\begin{equation}
\frac{d}{du}(e^{- u}\CN_\varepsilon(e^u)) = \frac{1}{2\pi} \int_\R \, e^{iut}C(1+it) \gamma_\varepsilon(t)dt\,.
\end{equation}
Compte tenu de (12), (15) donne
\begin{equation}
\frac{d}{du} (e^{-u} \CN_\varepsilon(e^u)) \in L^1(\R)\,.
\end{equation}
Donc $e^{-u}\CN_\varepsilon(e^u)$ a une limite quand $u\rg \infty$. De plus, comme
$$
\int \frac{d}{du}(e^{-u}\CN_\varepsilon(e^u)) du = C(1) \gamma_\varepsilon(0)\,,
$$
cette limite est $C(1)$ :
\begin{equation}
\lim_{x\rg \infty}
 \frac{\CN_\varepsilon(x)}{x} = C(1)\,.\end{equation}
Reste à passer de $\CN_\varepsilon(x)$ à $\CN(x)$.

Observons d'abord que
\begin{equation}
\CN_\varepsilon(e^u)= \int \CN(e^{u-v}) \hat{\gamma}_\varepsilon(v)dv
\end{equation}
avec pour $\hat{\gamma}
_\varepsilon(\cdot)$ la transformée de Fourier de $\gamma_\varepsilon(\cdot)$ :
$$
\hat{\gamma}_\varepsilon(v) = \frac{1}{\sqrt{2\pi}\varepsilon} \exp \big(-\frac{v^2}{2\varepsilon^2}\big)\,.
$$
Voici un lemme qui permet le passage.

\vspace{2mm}

\textbf{Lemme.} \textit{Si $\CM(x)$ est une fonction positive croissante sur $\R^+$, les propositions suivantes sont équivalentes :}

a) \textit{$\CM(x)$ a une densité (toujours sous--entendu : finie et non nulle)}

b) \textit{pour $\varepsilon>0$ assez petit, $\CM_\varepsilon(x)$ a une densité, $\CM_\varepsilon(e^u)$ étant la convolution de $\CM(e^u)$ et de $\hat{\gamma}_\varepsilon(u)$.}

\vspace{2mm}

\textbf{Preuve.} L'implication a) $\Rightarrow$ b) va de soi. Montrons que si a) n'a pas lieu, il en est de même pour b). Les cas $\overset{}{\underset{x\rg \infty}{\overline{\lim}}}\frac{\CM(x)}{x}=\infty$ et $\overset{}{\underset{x\rg \infty}{\lim}}\frac{\CM(x)}{x}=0$  sont évidents : ils entraînent respectivement $\overset{}{\underset{x\rg \infty}{\overline{\lim}}}\frac{\CM_\varepsilon(x)}{x}=\infty$ et $\overset{}{\underset{x\rg \infty}{\lim}}\frac{\CM_\varepsilon(x)}{x}=0$
 quel que soit $\varepsilon>0$.
Reste à étudier le cas
$$
\overset{}{\underset{x\rg \infty}{\underline\lim}}\frac{\CM(x)}{x} < \overset{}{\underset{x\rg \infty}{\overline{\lim}}}\frac{\CM(x)}{x} <\infty\,.
$$
Choisissons deux suites rapidement croissantes, $(x_n)$ et $(x_n')$ $(n\in \N)$ telles que
$$
\overset{}{\underset{n\rg \infty}{\lim}}\frac{\CM(x_n)}{x_n} = a < b = 
\overset{}{\underset{n\rg \infty}{\lim}}\frac{\CM(x_n')}{x_n'}\,.
$$
Posons $y_n=x_n(1-\delta)$ et $y_n'=x_n'(1+\delta)$ et fixons $a'$ et $b'$ tels que $a<a'<b'<b$. Pour $\delta$ assez petit et $n$ assez grand,
$$
\frac{\CM(y_n)}{y_n} <a' <b' < \frac{\CM(y_n')}{y_n'}\,.
$$
On voit alors en faisant la figure que, pour $\varepsilon$ assez petit,
$$
\overset{}{\underset{x\rg \infty}{\underline{\lim}}}\frac{\CM_\varepsilon(y_n)}{y_n}
< \overset{}{\underset{x\rg \infty}{\overline{\lim}}}\frac{\CM_\varepsilon(y_n')}{y_n'}
$$
donc b) n'a pas lieu.

Le lemme s'applique à $\CN(x)$ et $\CN_\varepsilon(x)$, et il résulte donc de (17) que $\CN(x)$ a une densité. D'après (17), cette densité est~$C(1)$.

On a donc établi le théorème de Diamond :

\vspace{2mm}

\textit{Sous la condition $(6)$, équivalente à $(2)$, la fonction $\CN(x)$ a une densité.}

\vspace{2mm}

\textbf{Extension et récapitulation}

L'algèbre de Wiener $W=\CF L^1(\R)$ est bien adaptée à la question de la densité de $\CN(x)$. Mais la méthode s'applique à d'autres algèbres de Banach de la forme $V=\CF H(\R)$. On conserve les notations $\CP(y)$, $\CN(x)$, $Z(s)$, $a(y)$ avec les hypothèses (4) et (6), $A(s)$, $B(s)$, $C(s)$ et $\CN_\varepsilon(x)$. On prend pour $H(\R)$
une algèbre de convolution de $R$ qui est 1) une algèbre de Banach 2) telle que
les fonctions de Gauss $\gamma_\varepsilon(\cdot)$ appartiennent à $H(\R)$ et que l'espace $\CD(\R)$ des fonctions indéfiniment dérivables à support compact s'injecte continument dans $H(\R)$ : cela garantit l'existence des partitions de l'unité qui permettent le passage du local au global (c'est une des idées maitresses de Wiener). En plus de (6) on suppose maintenant
\begin{equation}
e^{-u} a(e^u) \in H(\R)\,.
\end{equation}
On obtient de nouveau (12), (13), (14),(15) et enfin (16) avec $H(\R)$ au lieu de $L^1(\R)$, soit
\begin{equation}
\frac{d}{du} (e^{-u} \CN_\varepsilon(e^u)) \in H(\R)
\end{equation}
pour tout $\varepsilon>0$. Il n'en découle pas que $\frac{d}{du}(e^{-u}\CN(e^u)) \in H(\R)$ et  d'ailleurs c'est faux pour $H(\R)= L^1(\R)$ et $\CN=\N^+$. Le lemme assure seulement que $e^{-u}\CN(e^u)$ a une limite finie et non nulle quand $u\rg \infty$. Résumons, en deux parties

\vspace{2mm}

\textbf{Théorème 1.} a) \textit{L'hypothèse $(2)$ entraîne la conclusion $(1)$, avec $D=C(1)$ (théorème de Diamond). $b)$  Soit $H(\R)$ une algèbre de convolution sur $\R$  satisfaisant aux conditions $1)$ et $2)$ ci--dessus. On suppose $(19)$, c'est--à--dire}
\begin{equation}
e^{-u}(\CP(e^u) -\CP_o(e^u)) \in H(\R)\,.
\end{equation}
\textit{Alors, $\CN_\varepsilon(e^u)$ étant la convolution de $\CN(e^u)$ et de $\hat{\gamma}_\varepsilon(u)=\frac{1}{\sqrt{2\pi}} \exp(-\frac{u^2}{2\varepsilon^2})$, $e^{-u} \CN_\varepsilon(u)$ est une primitive d'une fonction appartenant à $H(\R)$, quel que soit~$\varepsilon>0$.}

\vspace{2mm}

La conclusion de la partie b) exprime une propriété globale peu élégante sur la fonction $\CN(e^u)$. La propriété locale qui en découle est plus parlante. En effet, si la déconvolution de $\CN_\varepsilon(e^u)$ par $\hat{\gamma}
_\varepsilon(u)$ conserve localement l'appartenance à $H(\R)$, (21) entraîne que, localement, $\CN(e^u)$ est primitive d'une fonction appartenant à $H(\R)$.

\section{Une autre condition, provenant aussi de l'analyse harmonique}

Nous abandonnons maintenant l'hypothèse (6), en supposant toujours (4), donc $a(y)=o(y)$ $(y\rg \infty)$. Les formules (7) à (10) définissent des fonctions  analytiques dans le demi--plan ouvert $\sigma>1$, et de nouveau $A(s)=O(|s|)$, $B(s)=O(|s|)$ et $C(s) = \exp O(|s|)$ quand $|s|\rg \infty$. Pour le prolongement de ces fonctions au demi--plan fermé $\sigma\ge 1$ on a besoin d'une nouvelle hypothèse sur $a(y)$. Rappelons~que
\begin{equation}
A(s) = \int y^{-s} da(y)\,.
\end{equation}
Notre hypothèse sera que la fonction $A(1+it)$ existe $(t \in \R)$, qu'elle est bornée et qu'elle vérifie la condition
\begin{equation}
|A(1+it) -A(1+it') | \le \omega (|t'-t|) \quad (t\in \R, \ t'\in \R)\,,
\end{equation}
la fonction $\omega(\cdot)$ étant croissante et telle que
\begin{equation}
\int_0^1 \frac{\omega(t)}{t} dt < \infty
\end{equation}
Alors la formule (22) définit $A(s)$ comme fonction continue et bornée sur le demi--plan fermé $\sigma\ge 1$. Il en est de même des fonctions $B(s)$ et $C(s)$, qui vont vérifier également des inégalités du type (23). Ainsi, il existe une constante $K$ telle que, pour $t\in \R$, $t'\in \R$ et $\sigma\ge 1$
\begin{equation}
|C(\sigma+it') -C(\sigma+it)\le K\, \omega(|t'-t|)\,.
\end{equation}

Or nous pouvons reprendre en la précisant la formule donnant $\CN(x)$ en fonction de $Z(s)$. Comme, pour~$\sigma>1$,
$$
\frac{Z(\sigma+it)}{\sigma{+it}} = \int e^{-iut} e^{-\sigma u} \CN(e^u) du\,,
$$
la formule d'inversion de Fourier donne formellement
$$
e^{-\sigma u} \CN(\varepsilon^u) = \frac{1}{2\pi} \int e^{iut} \frac{Z(\sigma+it)}{\sigma+it} dt \,;
$$
cela signifie qu'il y a un procédé de sommation du second membre qui donne le premier. Cette formule s'applique aux fonctions $\CN_o(\cdot)$ et $Z_o(\cdot)$ définies en (5), et on va utiliser le fait que $Z(s) =Z_o(s)C(s)$ (formule (10)). Ainsi, pour $\sigma>1$.
\begin{eqnarray}
\ \ e^{-\sigma u} (\CN(e^u)\! -\! C(\sigma) \CN_o (e^u)) \kern-3mm&=&\kern-3mm\displaystyle\frac{1}{2\pi}\int\! e^{iut} \frac{Z_o(\sigma+it)}{\sigma+it} (C(\sigma+it)\!-\!C(\sigma))\ dt\\
\kern-4mm&=&\kern-3mm \displaystyle \frac{1}{2\pi} \int e^{iut} \frac{C(\sigma+it)-C(\sigma)}{\sigma+it-1}\ dt\nonumber
\end{eqnarray}
D'après (25), les intégrales sont maintenant des intégrales de Lebesgue, et le module de l'intégrante dans (26) est majoré par $K$. On peut donc étendre la formule (26) à $\sigma=1$
\begin{equation}
e^{-u}(\CN(e^u)-C(1) \CN_o(e^u)) = \frac{1}{2\pi} \int e^{iut} \frac{C(1+it)-C(1)}{it}\ dt
\end{equation}
et
$$
\frac{C(1+it)-C(1)}{Ct} \in L^1(\R,dt)\,.
$$
Le second membre de (27) tend vers 0 quand $u\rg \infty$, et $e^{-u}\CN_o(e^u)=1$ pour $u>0$. Il en résulte que
\begin{equation}
\lim_{u\rg \infty} (e^{-u} \CN(e^u)) = C(1)\,.
\end{equation}

Enonçons le résultat.

\vspace{2mm}

\textbf{Théorème 2.} \textit{Supposons que $a(y)=o(y)$ $(y\rg \infty)$ et que la fonction $A(1+it)= \int y^{-1-it}da((y)$ existe, qu'elle est bornée et qu'elle vérifie $(23)$ et $(24)$. Alors $\CN(x)$ a une densité, à savoir~$C(1)$.}

\vspace{2mm}

Nous allons montrer que le théorème 2 peut s'appliquer dans des cas où l'hypothèse de Diamond n'est pas vérifiée.

\vspace{2mm}

\textbf{Exemples}

Choisissons
\begin{equation}
da = \Sigma \pm a_n \delta_{e^n}\ (n\in \N)
\end{equation}
avec des $a_n$ positifs qui seront choisis plus tard. Alors
$$
A(1+it) =\Sigma \pm a_n e^{-n} e^{-int}\,.
$$
Les signes $\pm$ étant choisis au hasard avec probabilité $\frac{1}{2}$ et indépendamment les uns des autres, on connait des conditions sur les $a_n$ pour que $A(1+it)$ existe, soit bornée, et vérifie (23) et (24) presque sûrement\footnote{Voir par exemple le chapitre 7 de \textit{Some random series of functions} by J.--P.~Kahane, Cambridge University Press~1985}. Il en est ainsi~si
$$
a_n e^{-n}= n^{-\alpha}
$$
avec $\frac{1}{2}<\alpha\le 1$ parce qu'alors $A(1+it)$ est p.s. hölderienne d'ordre $\beta$ pour tout $\beta<\alpha-1/2$. Désormais nous faisons ce choix. Alors le théorème 2 s'applique et, presque sûrement, $\CN(x)$ a une densité.

Calculons $a(y)$. En posant $h(\cdot)=1_{\R^+}(\cdot)$,
$$
a(y) = \Sigma \pm h(y-e^{n})a_n.
$$
Pour $y$ entre $e^n$ et $e^{n+1}$,
\begin{eqnarray*}
|a(y)| > a_n - \displaystyle \sum_{m<n} a_m> \frac{1}{10} e^n n^{-\alpha}\\
\displaystyle \int_{e^n}^{e^{n+1}}\frac{|a(y)|}{y^2}\ dy > \frac{1}{10}\ e^n n^{-\alpha} \int_{e^n}^{e^{n+1}} \frac{dy}{y^2}
\end{eqnarray*}
donc
$$
\int \frac{|a(y)|}{y^2}\ dy = \infty\,.
$$
Le théorème de Diamond ne s'applique pas et $\CN(x)$ a une densité.

Il y a beaucoup de variantes de (28). En voici une, facile à exprimer, qu'on obtient en changeant un peu la définition de $a(y)$, à savoir
$$
a(y) = \CP(y) -\CP_1(y)\,,\ \CP_1(y)=\pi(y)\,,
$$
et en choisissant alors dans (28) $a_n=\pi(e^n)$. $\CP$ est alors le multiensemble qu'on obtient à partir de $\CP_1$, l'ensemble des nombres premiers usuels, en décidant à pile ou face, pour chaque intervalle $[{e}^n,e^{n+1}[$, soit d'y supprimer tous les nombres premiers usuels qu'il contient, soit de les compter deux fois. On a la même conclusion que précédemment.

\section{Un cadre où la condition de Diamond est optimale}

Nous avons montré que la condition de Diamond (2) est suffisante et n'est pas nécessaire pour que la fonction $\CN(x)$ ait une densité (finie et non nulle). Cependant, dans des cas très particuliers, elle est nécessaire et suffisante : voir les observations 1 et 2. Peut--on  dans cette direction donner un énoncé général ? C'est ce que nous allons faire maintenant.

Rappelons qu'on associe à un ensemble $\CP$ multiplicativement libre (pour simplifier) porté par $]1,\infty[$ et localement fini, l'ensemble $\CN$ qu'il engendre par multiplication. En désignant par $\CP(y)$ et $\CN(x)$ les fonctions de comptage de $\CP$ et de $\CN$, on a la relation (3), que Beurling adopte pour définir $\CN(x)$ quand $\CP(y)$ désigne une fonction croissante arbitraire sur $]1,\infty[$. Diamond a montré que (2) entraîne (1), et on exprime (1) en disant que $\CN$, ou $\CN(x)$, a une densité.

L'exemple canonique sera ici $\CP=\CP_1$, l'ensemble des nombres premiers usuels, et $\CN=\CN_1=\N^+$, l'ensemble des entiers naturels $>0$. La condition de Diamond exprime une certaine proximité de $\CP$ à $\CP_1$ ou de $\CP(y)$ à $\CP_1(y)$ $(=\pi(x))$, à savoir
\begin{equation}
\int \frac{|\CP(y)-\CP_1(y)|}{y^2} \ dy < \infty\,.
\end{equation}
On peut modifier $\CP(y)$ en conservant (30). En effet :

\vspace{2mm}

\noindent\textbf{règle 1 :} si la condition est satisfaite pour $\CP(y)$, elle l'est pour toute fonction croissante de la forme $\CP(y)+O(1)$ $(y\rg \infty)$ ;

\noindent\textbf{règle 2 :} si la condition est satisfaite pour $\CP(y)$, elle l'est pour toute fonction croissante $\CP'(y)$ telle que $|\CP'(y) - \CP_1(y)| \le |\CP(y) - \CP_1(y)|\ (y>1)$.

\vspace{2mm}

\textbf{Théorème 3.} \textit{Parmi les conditions qui vérifient les règles $1$ et $2$ et qui entraînent que $\CN(x)$ a une densité (formule $(1)$), la condition de Diamond est optimale : toutes les autres l'impliquent.}

\vspace{2mm}

On peut dire également que si on impose les règles 1 et 2, la condition de Diamond est nécessaire et suffisante pour que $\CN(x)$ ait une densité. La suffisance étant établie, cela signifie ceci : si $\CP(y)$ est telle que pour toutes les fonctions $\CP'(y)$ qu'on en tire au moyen des règles 1 et 2, la fonction $\CN(x)$ a une densité, alors $\CP(y)$ vérifie~(30).

\vspace{2mm}

\textbf{Démonstration}

Supposons le contraire de (30), à savoir
\begin{equation}
\int \frac{|\CP(y)-\CP_1(y)|}{y^2}\ dy = \infty
\end{equation}
pour une fonction croissante $\CP(y)$ $(y>1)$ telle que $\CP(y)=o(y)$ $(y\rg \infty)$. Nous voulons montrer qu'en modifiant de façon convenable $\CP(y)$ suivant les règles 1 et 2, on parvient à une nouvelle fonction $\CP(y)$ telle que la fonction $\CN(x)$ associée ne vérifie pas (1). D'abord, suivant la règle 1, nous pouvons modifier la fonction $\CP(y)$ donnée de façon qu'elle devienne la fonction de comptage d'un ensemble $\CP$ pour lequel $\CP\cup\CP_1$ est multiplicativement libre, avec $\CP\cap \CP_1=\varnothing$. Cela étant, considérons les deux fonctions $(\CP(y)-\CP_1(y))^+$ et $(\CP(y)-\CP_1(y))^-$. Leur somme est $|\CP(y)-\CP_1(y)|$, donc l'une au moins peut être substituée à $|\CP(y)-\CP_1(y)|$ dans (30). Pour fixer les idées, disons que c'est $(\CP(y)-\CP_1(y))^+$. Soit $\CP^+(y)$ la fonction de comptage d'un ensemble $\CP^+$ qui coïncide avec $\CP$ là où $\CP(y)>\CP_1(y)$, et avec $\CP_1$ là où $\CP(y)<\CP_1(y)$. Ainsi
$$
(\CP(y)-\CP_1(y))^+ = \CP^+(y)- \CP_1(y)\,,
$$
et $\CP^+(y)$ s'obtient à partir de $\CP(y)$ en respectant la règle 2. Nous avons
$$
\int \frac{\CP^+(y)-\CP_1(y)}{y^2}\, d\, y = \infty
$$
et cela s'écrit aussi
\begin{equation}\label{F32}
\lim _{y\rg \infty} \Big(\sum_{p^+\in \CP^+\cap]1,Y]} \frac{1}{p^+} - \sum_{p\in \CP_1 \cap]1,Y]} \frac{1}{p}\Big) = \infty\,.
\end{equation}
Partons de $\CN_1$, l'ensemble des entiers naturels, engendré multiplicativement par $\CP_1$, et construisons l'ensemble $\CN^+$ associé à $\CP^+$ comme limite de $\CN^+[Y]$ ainsi définis : $\CN^+[Y]$ est associé à $\CP^+[Y]$, qui s'obtient à partir de $\CP_1$ en supprimant tous les $p\in \CP_1\cap]1,Y]$ et en ajoutant tous les $p^+\in \CP^+\cap ]1,Y]$. Or, vu l'indépendance multiplicative, chaque suppression d'un $p$ multiplie  la densité des entiers par $1-\frac{1}{p^+}$ et chaque ajout d'un $p^+$ la multiplie par $(1-\frac{1}{p^1})^{-1}$. La densité de $\CN^+[Y]$ est donc
$$
\prod_{p^+\in \CP^+\cap]1,Y]} \Big(1-\frac{1}{p^+}\Big)^{-1} \prod_{p\in \CP_1\cap]1,Y]}\Big(1- \frac{1}{p}\Big)
$$
et  (\ref{F32}) montre qu'elle tend vers l'infini quand $Y\rg \infty$. Donc
\begin{equation}
\lim _{x\rg \infty} \frac{\CN^+(x)}{x}=\infty\,.
\end{equation}
C'est bien la conclusion désirée, sous l'hypothèse qu'on peut substituer $(\CP(y) - \CP_1(y)^+$ à $|\CP(y)-\CP_1(y)|$ dans (30). Si ce n'est pas le cas, on travaille avec. $(\CP(y)-\CP_1(y))^-$, on définit $\CP^-(y)$ et on aboutit à
\begin{equation}
 \lim _{Y\rg\infty} \Big(\sum_{p\in \CP_1\cap]1,Y]}\frac{1}{p} - \sum_{p^-\in \CP^-\cap]1,Y]}\frac{1}{p^-}\Big)=\infty
 \end{equation} 
 au lieu de (32), puis à
 \begin{equation}
 \lim _{x\rg \infty} \frac{\CN^-(x)}{x} = 0
 \end{equation}
 au lieu de (33). Le théorème est donc démontré.
 
 \vspace{2mm}
 
\textbf{Remarque}

C'est pour la simplicité de l'argument que nous avons décomposé $|\CP(y)-\CP_1(y)|$ sous la forme $(\CP(y)-\CP_1(y))^+ + (\CP(y)-\CP_1(y))^-$ sur toute la demi--droite $]1,\infty[$, puis que nous avons défini $\CP^+(y)$ et $\CP^-(y)$ pour que $\CP^+(y)-\CP_1(y) =(\CP(y)-\CP_1(y))^+$ et $\CP_1(y)-\CP^-(y)=(\CP(y)-\CP_1(y))^-$. On pourrait définir des $\CP'(y)$ égaux à $\CP^+(y)$ sur certains grands intervalles, et à $\CP^-(y)$ sur des intervalles complémentaires. Au lieu de (33) et (35), on peut obtenir ainsi
$$
\overset{}{\underset{x\rg \infty}{\underline{\lim}}}
\frac{\CN'(x)}{x} < \overset{}{\underset{x\rg \infty}{\overline{\lim}}} \frac{\CN'(x)}{x}\,,
 $$
 ce qui est une autre façon de contredire (1).

\section{Sources et remerciements}

Cet article a pour origine l'usage qui est fait du théorème de Diamond dans \cite{4}. Ce théorème apparait dans \cite{1} au côté d'autres qui lient les propriétés des nombres premiers de Beurling et des entiers qu'ils engendrent. Michel Balazard distingue ces théorèmes selon le type de démonstration, analyse harmonique ou attaque directe qu'il appelle, suivant l'usage en théorie des nombres, méthode élémentaire. Ici le théorème de Diamond et le théorème~2 sont établis par analytique harmonique, et le théorème~3 par attaque directe.

L'article est signé de mon seul nom, et c'est une anomalie, parce qu'il fait suite à un travail en collaboration avec Eric Saïas, sa motivation et qu'il porte sa marque aussi bien pour se motivation que pour ses idées. Du moins est--ce l'occasion pour moi de dire ce que je dois à Eric au cours de ces dernières années, à cause de son enthousiasme, de sa curiosité, de ses connaissances et de son jugement. Merci, Eric

\vskip1cm

\vskip4mm
\begin{tabular}{p{6cm}l}

&Jean--Pierre Kahane  \\

&Laboratoire de Math\'ematiques d'Orsay \\

&Universit\'e Paris--Sud, CNRS \\

&Universit\'e Paris--Saclay \\

&91405 Orsay (France)  \\

\vspace{2mm}

&\textsf{jean-pierre.kahane@u-psud.fr} 
\end{tabular}

\end{document}